\documentstyle[12pt,leqno]{article}
\newtheorem{theorem}{\quad\large Theorem}

\begin{document}
\title{Three Sampling Formulas\\
}
\author{Alexander V. Gnedin\\ {\it Utrecht University}
}
\date{}
\maketitle

\par {\bf Abstract.} Sampling formulas describe probability laws of exchangeable combinatorial structures
like partitions and compositions.
We give a brief account of two known parametric families of sampling formulas and add a new
family to the list.

\vskip0.5cm

\par {\bf 1 Introduction.} By an integer composition of weight $n$ and length $\ell$ we shall mean
an ordered collection of 
positive integer
parts $\lambda =(\lambda_1,\ldots,\lambda_{\ell})$; we write $\lambda\vdash n$ for $\sum\lambda_j=n$.
It will be convenient to also use variables
$\Lambda_k=\lambda_k+\ldots +\lambda_{\ell},$ $k\leq \ell$,
so that $\lambda_j=\Lambda_j-\Lambda_{j-1}.$ 
\par A {\it composition structure} is a nonnegative function $q$ on compositions such that for each $n$
the values $\{q(\lambda):\lambda\vdash n\}$ comprise a probability distribution, say
 $q_n$, and the $q_n$'s satisfy
the following sampling consistency condition. Imagine an ordered series of randomly many 
nonempty boxes filled in randomly
with balls, so that the distribution of occupancy numbers from the left to the right is $q_n$.
The condition requires that if some $k<n$ balls are sampled out uniformly at random then the distribution
of the reduced occupancy numbers 
in nonempty boxes (in same order) must be exactly $q_{n-k}$ (without loss of generality we can take $k=1$).

\par Ignoring the order of boxes yields Kingman's partition structure \cite{RPS} 
(see \cite{Aldous} and \cite{CSP} for 
systematic development of the theory of partition structures, their relation to exchangeability 
and many references). But the relation cannot be uniquely inverted,
because for a given partition structure there are many ways to introduce the order in a consistent fashion.

\par Gnedin \cite{RCS} showed that all  composition structures can be uniquely 
represented by a  version of the 
Kingman's paintbox construction \cite{RPS}.
Let $U$ be a {\it paintbox}  -- a random open subset of $[0,1]$. With a paintbox we associate an ordered 
partition of $[0,1]$ comprised of the intervals of $U$ and of individual elements of the complement $U^c$,
with the order of blocks induced by the order on reals.
Suppose $n$ independent uniform random points are sampled from $[0,1]$ independently of $U$.
The sample points group somehow within the partition blocks and we  obtain
a random composition by writing 
the nonzero 
occupation numbers from the left to the right.
With probability one there is no tie among the sample points and the
consistency for various sample sizes follows from exchangeability in the sample.

\par From a topological viewpoint, the representation establishes a homeomorphism between the space 
of extreme composition structures
and the compact space of open subsets of $[0,1]$ (endowed with a weak topology), and also 
identifies the generic composition structure
with a unique mixture of extremes.
Thus 
already the set of extremes is intrinsically infinite-dimensional, 
not to say about the mixtures.
It is therefore a question of interest to find 
smaller parametric families which admit a reasonably simple
description.
\par In this note we discuss briefly three such families: one is an ordered 
modification, due to Donnelly and Joyce \cite{DoJo},
of the ubiquitous Ewens sampling formula (corresponding to
$(0,\theta)$-partition structure from the Ewens-Pitman two-parametric family \cite{PTRF}); 
another one, due to Pitman \cite{Bernoulli}, is 
an ordered (symmetric) version 
of the $(\alpha,0)$-partition structure, and 
the third composition structure is new. 
Despite the fact that the new composition structure is, in a sense, constructed from the
beta distributions, like the first two,
the corresponding  partition structure does not fit in the Ewens-Pitman family.
All three belong to a large infinite-dimensional family of  regenerative compositions introduced
and characterised in \cite{GP} and all three are  mixed, i.e.
generated by genuinely random paintboxes.

\vskip0.5cm

\par {\bf 2 Ordered ESF.} In their encyclopaedical exposition of the multivariate Ewens distribution 
Ewens and Tavar{\'e} presented the ordered version of ESF (see \cite{ESF}, Eqn. 41.6),
\begin{equation}\label{e}
e(\lambda)=\frac{\theta^{\ell}\,n!\,}{[\theta]_n}\prod_{j=1}^{\ell} \frac{1}{\Lambda_j}\,\qquad \theta>0\,,
\end{equation}
in connection with a size-biased permutation of the Ewens partition structure (here and forth, $[\,\,\,]$ is the 
Pochhammer factorial).
The special case $\theta=1$  is well known to combinatorialists as 
the distribution of cycle lengths in a uniform random permutation, provided
the cycles are ordered by increasing of their least elements.

\par Donnelly and Joyce \cite{DoJo} observed 
that the formula also defines a 
 composition structure, i.e. that 
(\ref{e}) determines a consistent sequence of random partitions taken together
 with an {\it intrinsic} ordering of classes,
based neither on the sizes of classes nor on labeling of
 `balls in boxes'. 
They argued that the  ordered structure is of some significance for  
biological applications,
and proved the following paintbox  representation of $e$. 
\par Let $Z_j$ be independent random variables with beta density
$${\rm d}\omega=\theta \,(1-z)^{\theta-1}{\rm d}z\,.$$
Let $U_{e}$ be the open set complementary to the stick-breaking sequence 
$$1- \prod_{j=1}^k (1-Z_j)\qquad k=1,2,\ldots$$
taken together with the endpoints of $[0,1]$. Then rephrasing Theorem 10 from \cite{DoJo} we have
\vskip0.5cm
\begin{theorem} The composition structure $e$ can be derived from the paintbox $U_e$.
\end{theorem}

\vskip0.5cm
\par The proof of this result given in \cite{DoJo} relied on the twin fact about weak convergence of the
size-biased permutation of ESF. Next is a direct argument which offers some more insight and exemplifies
the approach taken in this paper.

\vskip0.5cm

{\it Proof.} Introduce the binomial moments
\begin{eqnarray}\label{w}
w(n:m)&=&{n\choose m} \int_0^1 z^m (1-z)^{n-m}{\rm \,d}\omega(z)\\
&=&\theta\,{n\choose m}{\rm \, B} (m+1,n-m+\theta)\nonumber
 \end{eqnarray}
For $I$ the leftmost interval of $U$ (adjacent to $0$) the size of $I$ equals $Z_1$. 
Denoting $\hat{e}$ the composition structure derived from $U_{e}\,$ we aim to show that
$\hat{e}=e.$ 
\par The argument is based on two facts. Firstly,
suppose $n$ uniform points have been sampled from $[0,1]$ and
$I$ occured to contain $m$ sample points, then conditionally given $I$ the configuration of 
other $n-m$ points is as if it were a uniform sample from $I^c$. The second fact is that 
given $I$ the set $U_e\setminus I$ is a scaled distributional copy of $U_{e}$, as it is clear from the
definition of the paintbox via stick-breaking.

\par Composition
$(\lambda_1,\ldots,\lambda_{\ell})$ can only appear if the interval $J$,  defined to be the
leftmost of the intervals of $U_{e}$ discovered by the sample,  
contains exactly $\lambda_1$ sample points. 
The chance that $J$ coincides with $I$ is 
$w(n:m)$ and in this case the composition derived from the piece of $U_e$ to the 
right from $J$ must be $(\lambda_2,\ldots,\lambda_{\ell}).$
Otherwise $\lambda$ can appear only if $I$ contains no sample points and all $n$ group within
$U_e\setminus J$ in accord with $\lambda$. 

\par Combining these facts we get equation 
$$e(\lambda_1,\ldots,\lambda_{\ell})=w(n:\lambda_1)\,e(\lambda_2,\ldots,\lambda_{\ell})+w(n:0)\,
e(\lambda_1,\ldots,\lambda_{\ell})$$
leading to the recursion
\begin{eqnarray*}
e(\lambda_1,\ldots,\lambda_{\ell})=\frac{w(n:\lambda_1)}{1-w(n:0)} \,e(\lambda_2\ldots,\lambda_{\ell}).
\end{eqnarray*}
which is solved as  
\begin{eqnarray}\label{rec}
\hat{e} (\lambda)=\prod_{j=1}^{\ell} q(\Lambda_j:\lambda_j)
\end{eqnarray}
where 
\begin{eqnarray}\label{q}
q(n:m)&:=&\frac{w(n:m)}{1-w(n:0)}\\ 
&=&\frac{\theta}{n}\, \frac{n!}{(n-m)!}\,\frac{\,\,\,[\theta]_{n-m}}{[\theta]_n\,\,\,\,}\nonumber \,.
\end{eqnarray}
Cancelling common factors we arrive at (\ref{e}), thus $\hat{e}=e.$ $\Box$
\vskip0.5cm

\par There is a canonical correspondence between composition structures and probability distributions
of exchangeable compositions of an infinite set $\{\underline{1},\underline{2},\ldots\}$ (see \cite{RCS}). 
In terms of the paintbox representation the composition
derived from  $U$ is obtained by sampling infinitely many uniform points and
then assigning objects $\underline{i}$ and $\underline{j}$ to distinct classes if the closed interval
spanned on the $i$th and the $j$th sample points has a nonempty intersection with $U^c$. 
\par The infinite composition associated with $e$, call it ${\cal E}$, has a simply ordered
collection of blocks, and the law of large numbers says that the asymptotic frequencies of the
blocks (in a growing sample) coincide with the sizes of stick-breaking residuals, from the left to the right.
When we view $\cal E$ from the perspective of  restrictions ${\cal E}_n$ on finite sets
$\{\underline{1},\ldots,\underline{n}\}$'s, the collection of blocks {\it stabilises} (with probability one) 
in the sense that for any $k$ no {\it new} block appearing 
in ${\cal E}_{n'}\,$, for $n'>n\,$, will interlace with the collection of the first $k$ blocks represented in 
${\cal E}_n\,$, provided $n$ is sufficiently large (a zero-one law). Compositions
with this property were called `representable' in \cite{DoJo} and the class of such compositions 
generated by a general stick-breaking paintbox was characterised in \cite{GP}.

\vskip0.5cm

\par {\bf 3 PSF.} Pitman's composition structure is given by  Eqn. (30) in \cite{Bernoulli}:

\begin{equation}\label{PSF}
p(\lambda)=\frac{n!\,\alpha^{\ell}}{[\alpha]_n} \prod_{j=1}^{\ell} \frac{[1-\alpha]_{\lambda_j}}{\lambda_j!}\qquad 
0<\alpha<1\,.
\end{equation}
This sampling formula was derived from the following paintbox representation.
\vskip0.5cm
\begin{theorem} The paintbox $U_p$ for $p$ is the union of  excursion intervals
of the Bessel bridge of dimension $2-2\alpha$.
\end{theorem}
\vskip0.5cm
\par Equivalently, the complement $U_p^c$ is the set of zeroes of the Bessel bridge on $[0,1]$. 
The case $\alpha=1/2$ corresponds to the Brownian bridge. 
\par In fact, $p$ is a conditional version of another Pitman's composition structure $p'$ derived 
from
the set of zeroes of a Bessel process (which has final meander interval adjacent to the rightpoint of $[0,1]$).
Pitman obtained a formula for $p'$ akin to (\ref{PSF}) (see \cite{Bernoulli}, Eqn. (28))
using selfsimilarity of the Bessel process and distribution of the length of meander interval.
The relation between the structures is that 
$$p(\lambda)={\rm const}(n)\,p'(\lambda,1)\qquad \lambda\vdash n\,.$$  
\par Gnedin and Pitman \cite{GP} give a
 characterisation of $p$ related to the observation that this composition structure is also of the product form
(similar to (\ref{rec})) with
$$
q(n:m)=-\frac{{\alpha\choose m}{-\alpha\choose n-m}}{{-\alpha\choose n}}
$$

\par For $\ell$ fixed, $p$ is a symmetric function of the parts. This reflects in that $U_p$ is {\it symmetric},
that is has component intervals `in random order' (in  \cite{Aldous} the
open sets with this kind of invariance are called `exchangeable interval
partitions'). Summing $p(\lambda)$ over distinct permutations of parts yields a function on 
integer partitions which is the $(\alpha,0)$-partition structure from the Ewens-Pitman family.
It follows that $p$ could be obtained from the partition structure 
by permuting the parts in uniform random order (this is the general device allowing 
to derive symmetric composition structures and symmetric open sets 
from their unordered relatives \cite{JMVA}).

\par Blocks of the Pitman's composition $\cal P$ on $\{\underline{1}, \underline{2},\ldots,\}$
are ordered like the set of rational numbers and a such cannot be labeled by integers 
consistently with their intrinsic
order.
This happens each time  a composition has  infinitely many blocks (almost surely) and is symmetric.
A consequence is that the
infinite composition $\cal P$ has no definite first, second, etc  or the last block, in particular the 
first (hence $k$th) block in ${\cal P}_n$ does not stabilise as $n$ grows.

\vskip0.5cm
\par {\bf 4 A new sampling formula.} Here is a new composition structure 

\begin{equation}\label{g}
g(\lambda)=\frac{n!}{[\theta]_n}\prod_{j=1}^{\ell} \frac{1}{\lambda_j\, 
h_{\,\theta}(\Lambda_j)}\qquad \theta>0
\end{equation}
where
$$
h_{\,\theta}(n)=\sum_{k=1}^{n}\frac{1}{\theta+k-1}
$$
are the generalised harmonic numbers which coincide with the partial sums of the harmonic series when $\theta=1$. 

\par To explain the genesis of the formula consider stick-breaking with the general beta density
\begin{equation}\label{beta}
{\rm d}{\omega}(z)={\rm const}\cdot\, z^{\alpha-1}(1-z)^{\theta-1}\qquad \alpha,\theta>0\,.
\end{equation}
The resulting paintbox generates a composition structure given by the RHS of (\ref{rec}) with
\begin{equation}\label{qg}
q(n:m)={n\choose m} \frac{[\alpha]_m\,[\theta]_{n-m}}{[\alpha+\theta]_n-[\theta]_n}
\end{equation}
where (\ref{qg})  is obtained like (\ref{q}) from the binomial moments of the beta density (\ref{beta})
(to see this just follow the lines in the proof of Theorem 1).
\par For general $\alpha$ and $\beta$ the induced composition structure cannot be expressed by a simple
product formula, because the denominator has no good factorisation.
One notable exception is the ESF
appearing when $\alpha=1$.
Another exception is the case $\alpha=0$ giving rise to $g$; but this should 
be interpreted properly because measure $\omega$ becomes infinite.

\vskip0.5cm
\begin{theorem}
When $\alpha\downarrow 0$ the stick-breaking composition structure directed by the beta density {\rm (\ref{beta})}
converges to $g$.
\end{theorem}
{\it Proof.} Expansions in powers of $\alpha$ start with
$$[\alpha]_m= \alpha\,(m-1)!+\ldots \,,\qquad     [\alpha+\theta]_n-[\theta]_n=\alpha\, h_{\,\theta}(n)+\ldots$$
therefore when $\alpha$ approaches $0$ we get 
\begin{equation}\label{qq}
q(n:m)=\frac{n!}{(n-m)!\, } \frac{\,\,\,[\theta]_{n-m} }{  [\theta]_{n}\,\,\,\,} \,\frac{1}{m\,h_{\,\theta}(n)}\,.
\end{equation}
which yields $g$ as in (\ref{rec}). $\Box$

\vskip0.5cm

\par Distribution (\ref{qg}) underlying $g$ is especially simple for $\theta=1$ when it gives a weight 
proportional to $m^{-1}$ to each $m=1,\ldots,n$.

\par To determine the paintbox representation for $g$ we will extend the classical stick-breaking procedure
by embedding the process into continuous time and allowing infinitely many breaks within any time interval.
Note that defining a composition structure via the RHS of (\ref{rec}), through the
 binomial moments of some measure
$\omega$ and
$$q(n:m)=\frac{w(n:m)}{w(n:1)+\ldots+ w(n:n)}$$
we need not require that the measure $\omega$ be finite and do need to only impose the condition
$$\int_0^1 z\,{\rm d}\omega (z)<\infty$$
to have all binomial moments finite for $1\leq m\leq n<\infty$. 

\par In particular, our $g$ appears when we take improper density
\begin{equation}\label{sing}
{\rm d}\omega (z)= z^{-1} (1-z)^{\theta-1}\,{\rm d}z
\end{equation}
 (see \cite{GP} for more examples).  
For this $\omega$ consider a planar Poisson process (PPP) in the infinite strip $[0,\infty]\times [0,1]$
with ${\rm Lebesgue}\times \omega$ as intensity measure. The PPP has countably many atoms
$(\tau_j,\xi_j)$ 
(we adopt the conventional fake labeling of atoms which is not intended to say that
$\tau_j$ or $\xi_j$  is a definite random variable  for  particular $j$),  
and each location on the abscissa is a concentration point for the set of atoms.
Define a  pure-jump process with increasing cadlag paths 
$$S_t=1-\prod_{(\tau_j,\xi_j):\tau_j\leq t} (1- \xi_j)\,$$
where the product is over all PPP atoms to the left from $t$.
For any $t$ the product converges because $z\,\omega({\rm d}z)$ is a finite measure.
The process $(S_t)$ is a geometric subordinator: for $t'>t$ the ratio
$(1-S_{t'})/(1-S_t)$ is independent of the partial path on $[0,t]$ and has same distribution as 
$1-S_{t'-t}$.

\par (The reader feeling more comfort with breaking sticks from the right to the left 
should translate paintbox formulas using involution $z\leftrightarrow 1-z$ and also mirror the sampling formulas.)

\vskip0.5cm

\begin{theorem} The paintbox $U_g$ representing $g$ is the complement to the closure of the 
random set $\{S_t: t> 0\}$,  which is the range of the geometric subordinator.
\end{theorem}

{\it  Proof.} Fix $\lambda\vdash n$ and consider a uniform sample of size $n$. The composition $\lambda$ 
appears when 
for some $\tau_j$ the interval $[0,S_{\tau_j}]$ contains $m$ sample points grouped in one component interval
of  $U_{g}\cap [0,S_{\tau_j}]$ and the composition on the remaining $(n-m)$ sample points
is $(\lambda_2,\ldots,\lambda_{\ell})$. From the properties of uniform distribution and because
PPP is ruled by a product measure follows that the composition structure induced by $U_g$ is of the product form 
as in
(\ref{rec}) and we only need to justify the formula (\ref{qq}) for $q$ which is the distribution of
the first part of composition of $n$.

\par To that end, let $\pi (t)$ be the probability that some $m$ sample points group in one interval of
$U_{g}\cap [0,S_{t}]$ and denote $\epsilon_1,\ldots ,\epsilon_n$ the increasing order statistics
of uniform sample. Considering a small time interval $[0,{\rm d}t]$ it is not hard to  
see that $\pi$ satisfies the differential equation
$$\pi'=-a \,\pi +b\,,\qquad \pi (0)=0$$ 
with constant coefficients 
$$a=E\omega[\epsilon_1,1]=w(n:1)+\ldots +w(n:n)\,\,\,\,{\rm and }\,\,\,\,
b=E\omega[\epsilon_m,\epsilon_{m+1}]=w(n:m)$$
(with 1 in place of $\epsilon_{m+1}$ in case $m=n$) 
where $w(n:m)$'s are the binomial moments of (\ref{sing}). Solving the equation we obtain
$\phi(t)=(b/a)(1-e^{-at})\to b/a=q(n:m)$, 
as $t\to\infty$ whence $q(n:m)=b/a$ and this is (\ref{qq}).     $\Box$

\vskip0.5cm

\par The infinite composition $\cal G$ associated with $g$ has infinitely many blocks, and the set of blocks
is order isomorhic to the set of rational numbers. Unlike Pitman's $\cal P$ it is not symmetric, i.e. $g$
is sensible to permutation of parts $\lambda_j$ when $\ell>1$,
and the representing paintbox $U_g$ is not an `exchangeable interval partition'.
A combinatorialist might find natural to view $g$ as a function on Young diagrams $(\Lambda_1,\ldots,\Lambda_{\ell})$ with strictly
decreasing parts.

\par Ignoring the order in $\cal G$ yields a novel partition structure. For no $\theta$ belongs this partition structure
to the Ewens-Pitman two-parameter family, which 
had covered practically all explicit sampling formulas known to date. 
The distinction can be seen immediately by comparing 
the probability of one-class partition,  our $g(n)=q(n:n)$ given by (\ref{qg}) versus
the analogous quantity computed via Eqn. (16) in \cite{PTRF}
(the formulas do  not match for $n>4$ whatever the values of parameters). 

\vskip0.5cm

\par {\small Taking other integer values of $\alpha$ in (\ref{beta})  leads to formulas involving products
of stereotypic polynomial factors, e.g. for $\alpha=2$ we have
$$
g_2(\lambda)=
\frac{n!\,\theta^{\ell}\,(1+\theta)^{\ell}}{ [\theta]_n }\prod_{j=1}^{\ell}\frac{\lambda_j+1}{\Lambda_j+2\theta+1} \,\,.   
$$
The resulting infinite compositions  have simply ordered blocks and thus are more in line with $\cal E$.

\vskip0.5cm 
\noindent
{\it Department of Mathematics, Utrecht University\\
Postbus 80010, 3508 TA Utrecht, The Netherlands}\\
gnedin@math.uu.nl

\end{document}